\documentclass[12pt,british]{article}

\usepackage{babel}
\usepackage{amsmath,amsfonts,amssymb,amsthm}
\usepackage[utf8]{inputenc}
\usepackage[T1]{fontenc}
\usepackage{enumerate}
\usepackage{textcomp,url}
\usepackage{eucal}

\usepackage[noBBpl]{mathpazo}
\usepackage{mathabx}

\usepackage[matrix,arrow]{xy}

\makeatletter

\renewcommand{\p@enumii}{}

\def\@enum@{\list{\csname label\@enumctr\endcsname}%
           {\usecounter{\@enumctr}\def\makelabel##1{
\normalfont\ignorespaces\emph{{##1}~}}
\setlength{\labelsep}{3pt}
\setlength{\parsep}{0pt}
\setlength{\itemsep}{0pt}
\setlength{\leftmargin}{0pt}
\setlength{\labelwidth}{0pt}
\setlength{\listparindent}{\parindent}
\setlength{\itemsep}{0pt}
\setlength{\itemindent}{0pt}
\topsep=3pt plus 1pt minus 1 pt}}

\makeatother

\renewcommand{\epsilon}{\ensuremath{\varepsilon}}
\renewcommand{\phi}{\ensuremath{\varphi}}

\renewcommand{\to}{\ensuremath{\longrightarrow}}

\newcommand{\N}{\ensuremath{\mathbb N}}
\newcommand{\Z}{\ensuremath{\mathbb Z}}
\newcommand{\dt}{\ensuremath{\mathbb D}^{2}}
\newcommand{\St}[1][2]{\ensuremath{\mathbb S}^{#1}}
\newcommand{\FF}{\ensuremath{\mathbb F}}
\newcommand{\F}[1][n]{\ensuremath{\FF_{{#1}}}}
\newcommand{\rp}{\ensuremath{\mathbb{R}P^2}}

\newcommand{\sn}[1][n]{\ensuremath{S_{{#1}}}}

\DeclareRobustCommand*{\up}[1]{\textsuperscript{#1}}
\renewcommand{\th}{\ensuremath{\up{th}}}
\newcommand{\ft}[1][n]{\ensuremath{\Delta_{#1}}}

\renewcommand{\ker}[1]{\ensuremath{\operatorname{\text{Ker}}\left({#1}\right)}}

\newcommand{\aut}[1]{\ensuremath{\operatorname{\text{Aut}}\left({#1}\right)}}

\newcommand{\quat}[1][8]{\ensuremath{\mathcal{Q}_{#1}}}

\makeatletter
\def\@map#1#2[#3]{\mbox{$#1 \colon\thinspace #2 \to #3$}}
\def\map#1#2{\@ifnextchar [{\@map{#1}{#2}}{\@map{#1}{#2}[#2]}}
\makeatother

\newcommand{\brak}[1]{\ensuremath{\left\{ #1 \right\}}}
\newcommand{\ang}[1]{\ensuremath{\left\langle #1\right\rangle}}

\newcommand{\setr}[2]{\ensuremath{\brak{#1 \,\left\lvert \, #2 \right.}}}

\newtheoremstyle{theoremm}{}{}{\itshape}{}{\scshape}{.}{ }{}

\theoremstyle{theoremm}
\newtheorem{thm}{Theorem}
\newtheorem{lem}[thm]{Lemma}
\newtheorem{prop}[thm]{Proposition}
\newtheorem{cor}[thm]{Corollary}

\newtheoremstyle{remarkk}{}{}{}{}{\scshape}{.}{ }{}

\theoremstyle{remarkk}
\newtheorem{defn}[thm]{Definition}
\newtheorem{rem}[thm]{Remark}
\newtheorem{rems}[thm]{Remarks}

\newtheoremstyle{comment}{}{}{}{}{\bfseries}{:}{ }{}
\theoremstyle{comment}

\newcommand{\reth}[1]{Theorem~\protect\ref{th:#1}}
\newcommand{\relem}[1]{Lemma~\protect\ref{lem:#1}}
\newcommand{\repr}[1]{Proposition~\protect\ref{prop:#1}}
\newcommand{\reco}[1]{Corollary~\protect\ref{cor:#1}}
\newcommand{\resec}[1]{Section~\protect\ref{sec:#1}}
\newcommand{\rerem}[1]{Remark~\protect\ref{rem:#1}}

\newcommand{\req}[1]{equation~(\protect\ref{eq:#1})}

\begin{document}

\title{Classification of the virtually cyclic subgroups of the pure braid groups of the projective
plane}
\author{DACIBERG~LIMA~GON\c{C}ALVES\\
Departamento de Matem\'atica - IME-USP,\\
Caixa Postal~66281~-~Ag.~Cidade de S\~ao Paulo,\\ 
CEP:~05311-970 - S\~ao Paulo - SP - Brazil.\\
e-mail:~\texttt{dlgoncal@ime.usp.br}\vspace*{4mm}\\
JOHN~GUASCHI\\
Laboratoire de Math\'ematiques Nicolas Oresme UMR CNRS~\textup{6139},\\
Universit\'e de Caen BP 5186,\\
14032 Caen Cedex, France.\\
e-mail:~\texttt{guaschi@math.unicaen.fr}}

\date{\today}

\begingroup
\renewcommand{\thefootnote}{}
\endgroup 


\maketitle

\begin{abstract}\noindent
\emph{We classify the (finite and infinite) virtually cyclic subgroups of the
pure braid groups $P_{n}(\rp)$ of the projective plane. The maximal finite subgroups of
$P_{n}(\rp)$ are isomorphic to the quaternion group of order $8$ if $n=3$, and to $\Z_{4}$ if
$n\geq 4$. Further, for all $n\geq 3$, up to isomorphism, the following groups are the infinite
virtually cyclic subgroups of $P_{n}(\rp)$: $\Z$, $\Z_{2} \times \Z$ and the amalgamated product
$\Z_{4} \ast_{\Z_{2}} \Z_{4}$.}
\end{abstract}

\section{Introduction}\label{sec:intro}
 
The braid groups $B_n$ of the plane were introduced by E.~Artin
in~1925~\cite{A1,A2}. Braid groups of surfaces were studied by
Zariski~\cite{Z}. They were later generalised by Fox to braid
groups of arbitrary topological spaces via the following
definition~\cite{FoN}. Let $M$ be a compact, connected surface, and
let $n\in\N$. We denote the set of all ordered $n$-tuples of distinct
points of $M$, known as the \emph{$n\th$ configuration space of $M$},
by:
\begin{equation*}
F_n(M)=\setr{(p_1,\ldots,p_n)}{\text{$p_i\in M$ and $p_i\neq p_j$ if $i\neq j$}}.
\end{equation*}
Configuration spaces play an important r\^ole in several branches of mathematics and have been
extensively studied, see~\cite{CG,FH} for example. 

The symmetric group $\sn$ on $n$ letters acts freely on $F_n(M)$ by
permuting coordinates. The corresponding quotient will be denoted by
$D_n(M)$. The \emph{$n\th$ pure braid group $P_n(M)$} (respectively
the \emph{$n\th$ braid group $B_n(M)$}) is defined to be the
fundamental group of $F_n(M)$ (respectively of $D_n(M)$). 

Together with the $2$-sphere $\St$, the braid groups of the
real projective plane $\rp$ are of particular interest, notably because they have
non-trivial centre~\cite{vB,GG2}, and torsion elements~\cite{vB,M}.
Indeed, Van Buskirk showed that among the braid groups of compact,
connected surfaces, $B_n(\St)$ and $B_n(\rp)$ are the only ones to
have torsion~\cite{vB}. Let us recall briefly some of the properties
of $B_n(\rp)$~\cite{GG2,M,vB}.

If $\dt\subseteq \rp$ is a topological disc, there is a group homomorphism $\map
{\iota}{B_n}[B_n(\rp)]$ induced by the inclusion. If $\beta\in B_n$ then we shall denote its image
$\iota(\beta)$ simply by $\beta$. A presentation of $B_{n}(\rp)$ was given in~\cite{vB}, and of
$P_{n}(\rp)$ in~\cite{GG8}. The first two braid groups of $\rp$ are finite: $B_1(\rp)=P_1(\rp)\cong
\Z_{2}$, $P_{2}(\rp)$ is isomorphic to the quaternion group $\quat$ of order~$8$ and $B_{2}(\rp)$ is
isomorphic to the generalised quaternion group of order $16$. For $n\geq 3$, $B_{n}(\rp)$ is
infinite. The pure braid group $P_{3}(\rp)$ is isomorphic to a semi-direct product of a free group of
rank $2$  by $\quat$~\cite{vB}; an explicit action was given in~\cite{GG2} (see also the proof of
\repr{p3noz4z}).

The so-called `full twist' braid of $B_n(\rp)$ is defined by $\ft= (\sigma_1\cdots\sigma_{n-1})^n$.
For $n\geq 2$, $\ft[n]$ is the unique element of $B_n(\rp)$ of order $2$~\cite{GG2}, and it
generates the centre of $B_n(\rp)$~\cite{M}. The finite order elements of $B_{n}(\rp)$ were
characterised by Murasugi~\cite{M} (see \reth{murasugi}, \resec{murasugi}), however their orders are
not clear, even for elements of $P_{n}(\rp)$. In~\cite{GG2}, we proved that for $n\geq 2$, the
torsion of $P_{n}(\rp)$ is $2$ and $4$, and that of $B_{n}(\rp)$ is equal to the divisors of $4n$ and
$4(n-1)$. 

The classification of the finite subgroups of $B_{n}(\St)$ and $B_{n}(\rp)$ is an interesting
problem, and helps us to better understand their group structure. In the case of $\St$, this was
undertaken in~\cite{GG5,GG9}. It is natural to ask which finite groups are realised as
subgroups of $B_n(\rp)$. As for $B_{n}(\St)$, one common property of such subgroups is
that they are finite periodic groups of cohomological period $2$ or $4$. Indeed, by Proposition~6 
of~\cite{GG2}, the universal covering $X$ of $F_n(\rp)$ is a finite-dimensional complex which has the
homotopy type of $\St[3]$. Thus any finite subgroup of $B_n(\rp)$ acts freely on $X$, and so has
period $2$ or $4$ by Proposition~10.2, Section~10, Chapter~VII of~\cite{Br}. Since $\ft$ is the
unique element of order~$2$ of $B_n(\rp)$, and it generates the centre $Z(B_n(\rp))$, the Milnor
property must be satisfied for any finite subgroup of $B_n(\rp)$. 

In \resec{finitepn}, we start by determining the maximal finite subgroups of $P_{n}(\rp)$:
\begin{prop}\label{prop:finitepn}
Up to isomorphism, the maximal finite subgroups of $P_{n}(\rp)$ are:
\begin{enumerate}[(a)]
\item $\Z_{2}$ if $n=1$.
\item $\quat$ if $n=2,3$.
\item $\Z_{4}$ if $n\geq 4$.
\end{enumerate}
\end{prop}
In \resec{murasugi}, we simplify Murasugi's characterisation of the torsion elements of
$B_{n}(\rp)$ (see \repr{power}), which enables us to show that within $B_{n}(\rp)$, there are two
conjugacy classes of subgroups isomorphic to $\Z_{4}$ lying in $P_{n}(\rp)$ (see \repr{order4}).

The rest of the paper is devoted to determining the infinite virtually cyclic subgroups of
$P_{n}(\rp)$ (recall that a group is said to be \emph{virtually cyclic} if it has a cyclic subgroup
of finite index). This was initially motivated by a question of Stratos Prassidis concerning the
determination of the algebraic $K$-theory of the braid groups of $\St$ and $\rp$. It has been shown
recently that the full and pure braid groups of these surfaces satisfy the Fibered Isomorphism
Conjecture of T.~Farrell and L.~Jones~\cite{BJPL,JPM1,JPM2}. This implies that the algebraic
$K$-theory groups of their group rings may be computed by means of the algebraic $K$-theory groups of
their virtually cyclic subgroups via the assembly maps. More information on these topics may be found
in~\cite{BLR, FJ,JP}. As well as helping us to better understand these braid groups, this provides us
with additional reasons to find their virtually cyclic subgroups. 

In \resec{virtcyc}, we recall the criterion due to Wall of an infinite virtually cyclic group $G$ as
one which has a finite normal subgroup $F$ such that $G/F$ is isomorphic to $\Z$ or to the free
product $\Z_{2} \ast \Z_{2}$ (this being the case, we shall say that $G$ is of \textit{Type~I} or
\textit{Type~II} respectively). Wall's result enables us to establish a list of the possible
infinite virtually cyclic subgroups of a given group $G$, if one knows its finite subgroups
(which is the case for $P_{n}(\rp)$). The real difficulty lies in deciding whether the groups
belonging to this list are effectively realised as subgroups of $G$. In \relem{type2crit} we give a
useful criterion to decide whether a group is of Type~II. 

In the case of $P_{n}(\St)$, since there are only two finite subgroups for $n\geq 3$, the trivial
group and that generated by the full twist (which are the elements of the centre of $P_{n}(\St)$),
it is then easy to see that its infinite virtually cyclic subgroups are isomorphic to $\Z$ or
to $\Z_{2}\times \Z$.

Since the structures of the finite subgroups of $P_{3}(\rp)$ and $P_{n}(\rp)$ differ, we separate the
discussion of their virtually cyclic subgroups. However, it turns out that up to isomorphism, they
have the same infinite virtually cyclic subgroups:
\begin{thm}\label{th:vcp}
Let $n\geq 3$. Up to isomorphism, the infinite virtually cyclic subgroups of $P_{n}(\rp)$ are $\Z$,
$\Z_{2}\times \Z$ and $\Z_{4} \ast_{\Z_{2}} \Z_{4}$.
\end{thm}
For $n= 3$, the result will be proved in \reth{vcp3} (see \resec{virtcycp3}), while for $n\geq 4$, it
will be proved in \reth{vcp4} (see \resec{virtcycp4}).

As a immediate corollary of \repr{finitepn} and \reth{vcp}, we obtain the classification of
the virtually cyclic subgroups of $P_{n}(\rp)$:
\begin{cor}
Let $n\in \N$. Up to isomorphism, the virtually cyclic subgroups of $P_{n}(\rp)$ are:
\begin{enumerate}[(a)]
\item The trivial group $\brak{e}$ and $\Z_{2}$ if $n=1$.
\item $\brak{e}$, $\Z_{2}$, $\Z_{4}$ and $\quat$ if $n=2$.
\item $\brak{e}$, $\Z_{2}$, $\Z_{4}$, $\quat$, $\Z$,
$\Z_{2}\times \Z$ and $\Z_{4} \ast_{\Z_{2}} \Z_{4}$ if $n=3$.
\item $\brak{e}$, $\Z_{2}$, $\Z_{4}$, $\Z$,
$\Z_{2}\times \Z$ and $\Z_{4} \ast_{\Z_{2}} \Z_{4}$ if $n\geq 4$.
\end{enumerate}
\end{cor}

One of the key results needed in the proof of \reth{vcp}  is that $P_{n}(\rp)$ has no subgroup
isomorphic to $\Z_{4} \times \Z$. This fact allows us to eliminate several possible Type~I 
and Type~II subgroups, and has the following interesting corollary which we prove at the end of
\resec{virtcycp4}:
\begin{cor}\label{cor:centralx}
Let $n\geq 2$, and let $x \in P_{n}(\rp)$ be an element of order $4$. Then its centraliser
$Z_{P_{n}(\rp)}(x)$ in $P_{n}(\rp)$ is equal to $\ang{x}$.
\end{cor}

The study of the finite subgroups of $B_{n}(\rp)$ and of the infinite virtually cyclic subgroups of
$B_{n}(\St)$ and $B_{n}(\rp)$ is the subject of work in progress.

\subsection*{Acknowledgements}

This work took place during the visits of the second author to the Departmento de Matem\'atica do IME-Universidade de
S\~ao Paulo during the periods 10\up{th}~--~27\up{th}~August 2007 and 4\up{th}~--~23\up{rd}~October 2007,  and of the
visit of the first author to the Laboratoire de Math\'ematiques Emile Picard, Université Paul Sabatier, Toulouse, during
the period 27\up{th}~August~--~17\up{th}~September 2007. It was supported by the international Cooperation USP/Cofecub
project number 105/06, by an `Aide Ponctuelle de Coop\'eration' from the Universit\'e Paul Sabatier, by the
Pr\'o-Reitoria de Pesquisa  Projeto~1, and by the Projeto Tem\'atico FAPESP no.~2004/10229-6 `Topologia Alg\'ebrica,
Geometrica e Diferencial'. The authors would like to thank Stratos Prassidis for having brought the question of the
virtually cyclic subgroups to their attention.

\section{Finite subgroups of $P_{n}(\rp)$}\label{sec:finitepn}

In this section, we characterise the finite subgroups of $P_{n}(\rp)$ by proving \repr{finitepn}. In
Theorem~4, Corollary~19 and Proposition~23 respectively of~\cite{GG2}, we obtained the following
results:
\begin{thm}[\cite{GG2}]\label{th:agtresults}
Let $n\geq 2$. Then:
\begin{enumerate}[(a)]
\item $B_{n}(\rp)$ has an element of order $\ell$ if and only if $\ell$ divides either $4n$ or
$4(n-1)$.
\item\label{it:torpn} the (non-trivial) torsion of $P_{n}(\rp)$ is precisely $2$ and $4$.
\item the full twist $\ft$ is the unique element of $B_{n}(\rp)$ of order $2$.
\end{enumerate}
\end{thm}
It was shown in~\cite{M} that $\ft$ generates the centre of $B_{n}(\rp)$. Using the projection
$P_{n+1}(\rp)\to P_{n}(\rp)$ given by forgetting the last string, and induction, one may check that
$\ft$ also generates the centre of $P_{n}(\rp)$.

It follows from \reth{agtresults} that the maximal (finite) cyclic subgroups of $B_{n}(\rp)$ are
isomorphic to $\Z_{4n}$ or $\Z_{4(n-1)}$, and that those of $P_{n}(\rp)$ are isomorphic to $\Z_{4}$.

\begin{proof}[Proof of \repr{finitepn}.]
Since $P_{1}(\rp)\cong \Z_{2}$ and $P_{2}(\rp)\cong \quat$, the result follows easily for $n=1,2$.
So suppose that $n\geq 3$. Consider the Fadell-Neuwirth short exact sequence:
\begin{equation*}
1 \to P_{n-2}(\rp \setminus \brak{x_{1},x_{2}}) \to P_{n}(\rp) \stackrel{p_{\ast}}{\to} P_{2}(\rp) \to
1,
\end{equation*}
where $p_{\ast}$ corresponds geometrically to forgetting the last $n-2$ strings. Let $H$ be a finite
subgroup of $P_{n}(\rp)$. Then $p_{\ast}\lvert_{H}$ is injective. To see this, suppose that $x,y\in
H$ are such that $p_{\ast}(x)= p_{\ast}(y)$. Then $xy^{-1}\in H \cap \ker{p_{\ast}}$. But
$\ker{p_{\ast}}$ is torsion free, so $x=y$, which proves the claim. In particular, $\lvert H
\rvert \leq 8$, and since the torsion of $P_{n}(\rp)$ is $2$ and $4$, and $P_{n}(\rp)$ has a unique
element of order $2$, it follows that $H$ is isomorphic to one of $\brak{e}$, $\Z_{2}$ (so is
$\ang{\ft}$), $\Z_{4}$ or $\quat$. If $n=3$ then the above short exact sequence splits~\cite{vB}, 
and so $P_{3}(\rp) \cong \F[2] \rtimes \quat$. Hence the finite subgroups of $P_{3}(\rp)$ are those
of $\quat$. Now let $n\geq 4$. Then $\quat$ is not realised as a subgroup of $P_{n}(\rp)$. For
suppose that $H< P_{n}(\rp)$, where $H \cong \quat$. From above, it follows that
$p_{\ast}\lvert_{H}$ is an isomorphism, so $p_{\ast}$ admits a section. But this would contradict
Theorem~3 of~\cite{GG2}. The result follows from \reth{agtresults}(\ref{it:torpn}).
\end{proof}

\section{Murasugi's characterisation of the finite order elements of
$B_{n}(\rp)$}\label{sec:murasugi}

In this section, our aim is to reorganise the classification of Murasugi~\cite{M} of the
finite order elements of $B_{n}(\rp)$ in a form that shall be more convenient for our purposes. As
well as being useful in its own right, we shall make use of this reorganisation to
prove that there are precisely two conjugacy classes in $B_{n}(\rp)$ of subgroups of pure braids of
order $4$. This shall be applied in the classification of the virtually cyclic subgroups of
$P_{n}(\rp)$, notably to show that $P_{n}(\rp)$ has no subgroup isomorphic to $\Z_{4} \times \Z$ (see
Propositions~\ref{prop:p3noz4z} and~\ref{prop:noz4z}).

Up to conjugacy, the finite order elements of $B_{n}(\rp)$ may be characterised as
follows:
\begin{thm}[\cite{M}]\label{th:murasugi}
Every element of $B_{n}(\rp)$ of finite order is conjugate to one of the following elements:
\begin{align*}
A_{1}(n,r,s,q) &= (\rho_{r}\sigma_{r-1}\cdots \sigma_{1})^s (\sigma_{r+1}\cdots \sigma_{n-1})^q\\
A_{2}(n,r,s,q)&= (\rho_{r}\sigma_{r-1}\cdots \sigma_{1})^s (\sigma_{r+1}\cdots
\sigma_{n-1}\sigma_{r+1})^q,
\end{align*}
where $(n-r,q) \approx (2r,s)$ in the first case, and $(n-r-1,q) \approx (2r,s)$ in the second. The
relation $\approx$ is defined by:
\begin{equation*}
\text{$(a,b) \approx (c,d)$ if there exists $(m,k)\neq (0,0)$ such that $m(a,b)=k(c,d)$.}
\end{equation*}
\end{thm}

This characterisation may be simplified as follows:
\begin{prop}\label{prop:power}
Let $n\geq 2$. Then every finite order element of $B_{n}(\rp)$ is conjugate to a power of some
$A_{i}(n,r,2r/l,p/l)$, $i=1,2$, where 
\begin{equation*}
p=
\begin{cases}
n-r & \text{if $i=1$}\\
n-r-1 & \text{if $i=2$,}
\end{cases}
\end{equation*}
and $l=\gcd(p,2r)$.
\end{prop}

\begin{rem}
The characterisation of the finite order elements of $B_{n}(\rp)$ thus appears more transparent
than that given in \reth{murasugi}, in the sense that these elements may be determined directly from
$n$ and $r$ (without involving the integers $s,q$, nor the relation $\approx$).
\end{rem}

\begin{proof}[Proof of~\repr{power}]
From \reth{murasugi}, every finite order element of $B_{n}(\rp)$ is conjugate to some
$A_{i}(n,r,s,q)$, where $m(p,q) = k(2r,s)$. Thus every finite order element of $B_{n}(\rp)$ is
conjugate to a power of some $A_{i}(n,r,s/\lambda,q/\lambda)$, where $\lambda=\gcd(s,q)$. We
consider the following cases separately:
\begin{enumerate}[(a)]
\item If $r=0$ then $m=s=0$ and $k\neq 0$. Then $A_{i}(n,r,s/\lambda,q/\lambda)= A_{i}(n,0,0,1)= 
A_{i}(n,r,2r/l,p/l)$ as required.
\item If either $i=2$ and $r=n-1$, or $i=1$ and $r=n$ then $p=0$. Thus $k=0$ and $m\neq
0$, so $q=0$. Hence $A_{i}(n,r,s/\lambda,q/\lambda)= A_{i}(n,r,1,0)= A_{i}(n,r,2r/l,p/l)$ as
required.
\item Suppose that either $i\in \brak{1,2}$ and $1\leq r\leq n-2$, or $i=1$ and $r=n-1$, so $p\neq 0$.
Since $mp=2rk$, $mq=ks$ and $(m,k)\neq (0,0)$, it follows that $m,k\neq 0$, $m/k=s/q=2r/p$ (if
$s=0$ then $q=0$, but then $A_{i}(n,r,s,q)$ would be trivial), and so $s'p=2rq'$, where
$s'=s/\lambda$ and $q'=q/\lambda$ are coprime. We thus obtain $2r=ls'$ and $p=lq'$,
so $\frac{(2r,p)}{l}= \frac{(s,q)}{\lambda}$. Hence every finite order element of $B_{n}(\rp)$ is
conjugate to a power of some $A_{i}(n,r,2r/l,p/l)$.
\end{enumerate}
\end{proof}

In Proposition~26 of~\cite{GG2}, we proved that the following elements of $B_{n}(\rp)$:
\begin{align*}
a &= \sigma_{n-1}^{-1} \cdots \sigma_{1}^{-1} \rho_{1}\\
b &= \sigma_{n-2}^{-1} \cdots \sigma_{1}^{-1} \rho_{1}
\end{align*}
are of order $4n$ and $4(n-1)$ respectively. By Remark~27 of~\cite{GG2}, we have
\begin{equation}
\label{eq:defalpha}
\left\{ \begin{aligned}
\alpha &= a^n= \rho_{n} \cdots \rho_{1}\\
\beta &= b^{n-1}= \rho_{n-1} \cdots \rho_{1}.
\end{aligned} \right.
\end{equation}
It is clear that $\alpha$ and $\beta$ are pure braids of order $4$. Using \repr{power}, we now show
that any element in $P_{n}(\rp)$ of order $4$ is conjugate in $B_{n}(\rp)$ to one of these two
elements (or their inverses):
\begin{prop}\label{prop:order4}
Let $n\geq 2$. Then every element of order $4$ of $P_{n}(\rp)$ is conjugate in $B_{n}(\rp)$ to one
of $\alpha$ and $\beta$ or their inverses.
\end{prop}

\begin{rem}
Since the Abelianisation of $B_{n}(\rp)$ is isomorphic to $\Z_{2} \times \Z_{2}$ (cf.\ the
presentation of $B_{n}(\rp)$ given in~\cite{vB}), and the generator $\sigma_{i}$, $1\leq i \leq n-1$
(resp.\ $\rho_{j}$, $1\leq j \leq n$) is sent to the generator of the first (resp.\ second) $\Z_{2}$-factor, $\alpha$ and $\beta$ are not conjugate in
$B_{n}(\rp)$.
\end{rem}

\begin{proof}[Proof of \repr{order4}.]
For given $n$ and $r$, define
\begin{equation*}
\omega= 
\begin{cases}
\sigma_{r+1} \cdots \sigma_{n-1} & \text{if $i=1$}\\
\sigma_{r+1} \cdots \sigma_{n-1} \sigma_{r+1} & \text{if $i=2$.}
\end{cases}
\end{equation*}

Up to conjugacy and powers, the finite order elements of $B_{n}(\rp)$ appear in one of the three
cases given in the proof of \repr{power}. Among these elements, we search for pure braids of order
$4$. 
\begin{enumerate}[(a)]
\item We have that $A_{i}(n,0,0,1)= \omega$.
If $i=1$ then $\omega= \sigma_{1} \cdots \sigma_{n-1}$ is of order $2n$, and its permutation is
$(1,n,\cdots,2)$. The smallest $j\geq 1$ for which $\omega^j$ becomes pure is thus $j=n$, but then
$\omega^n=\ft$. If $i=2$ then $\omega= \sigma_{1} \cdots \sigma_{n-1}\sigma_{1}$ is of order
$2(n-1)$, its permutation is $(1,n,\cdots,3)$, and $\omega^{n-1}=\ft$. The smallest $j\geq 1$ for
which $\omega^j$ becomes pure is indeed $j=n-1$. In either case, the powers of $\omega$ yield no pure
braids of order $4$.

\item Suppose that either $i=2$ and $r=n-1$, or $i=1$ and $r=n$. Using the relation
$\rho_{i+1}=\sigma_{i}^{-1} \rho_{i} \sigma_{i}^{-1}$, $1\leq i\leq n-1$ in $B_{n}(\rp)$
(see~\cite{vB,M,GG2}), we have that
\begin{equation*}
A_{i}(n,r,1,0)= \rho_{r} \sigma_{r-1}\cdots \sigma_{1}=
\begin{cases}
b & \text{if $r=n-1$}\\
a & \text{if $r=n$.}
\end{cases}
\end{equation*}
Studying the permutations of $a$ and $b$, we observe that the smallest power of each of
these elements which is a pure braid is $b^{n-1}=\beta$ and $a^n=\alpha$ (which are of order $4$).

\item If $i=1$ and $r=n-1$ then we have $A_1(n,r,2r/l,p/l)= A_1(n,n-1,2(n-1),1)= (\rho_{n-1}
\sigma_{n-2}\cdots \sigma_{1})^{2(n-1)}=b^{2(n-1)}=\ft$ which is of order $2$. So we may
suppose that $1\leq r\leq n-2$. Let $\xi= \rho_{r} \sigma_{r-1}\cdots \sigma_{1}$ and
$y_{i}=A_{i}(n,r,2r/l,p/l)=\xi^{2r/l} \omega^{p/l}$. By equation~(5.1) at the bottom of page~79
of~\cite{M}, $y_{i}^l= \xi^{2r} \omega^{p}=\ft$, so the order of $y_{i}$ divides $2l$. Now the
permutation of $\omega$ is a $p$-cycle, and since $l$ divides $p$, the permutation of $\omega^{p/l}$
is a product of $p/l$ disjoint $l$-cycles. Thus the smallest $j\geq 1$ for which $(\omega^{p/l})^j$
becomes a pure braid is $j=l$. But $y_{i}^l$ is the full twist, and so the only pure braids among the
powers of $y_{i}$ are $\ft$ and the identity. In particular, no power of $y_{i}$ can yield a pure
braid of order $4$.
\end{enumerate}
We conclude that up to conjugacy and inverses, the only elements of $P_{n}(\rp)$ of order $4$ are
$\alpha$ and $\beta$, and that they only occur in the second case ($i=2$ and $r=n-1$, or $i=1$ and
$r=n$). This completes the proof of the proposition.
\end{proof}

As a corollary of the proof of \repr{order4}, we are able to determine explicitly the orders of the
torsion elements given by \repr{power} as follows.
\begin{cor}\label{cor:orderl}
Let $n\geq 2$. Then for all $i=1,2$, and all $r=0,1,\ldots,n$, the element
$y_{i}=A_{i}(n,r,2r/l,p/l)=\xi^{2r/l} \omega^{p/l}$ of $B_{n}(\rp)$ is of order $2l=2 \gcd(p,2r)$.
\end{cor}

\begin{proof}
If $1\leq r\leq n-2$, we deduce from the above proof that $y_{i}$ is of order $2l$. If $r=0$ then
$l=p$, and $y_{i}=\omega$ which is of order $2l$. If $i=1$ and $r=n-1$ then $l=1$, and $y_{i}=\ft$
is indeed of order $2$. Finally, if $i=2$ and $r=n-1$ or $i=1$ and $r=n$ then $p=0$ and $l=2r$.
Hence $y_{1}=\xi=a$ which is of order $4n$, and $y_{2}=\xi=b$ which is of order $4(n-1)$. The
result follows.
\end{proof}

Set $l_{1}(n,r)=\gcd(2r,n-r)$ and $l_{2}(n,r)=\gcd(2r,n-r-1)$. For $0\leq r\leq n-1$, clearly
$l_{2}(n,r)= l_{1}(n-1,r)$, so it suffices to know the values of $l_{1}(n,r)$. An alternative manner
to express the order $2l_{1}(n,r)$ of $y_{1}$ is given by the following:
\begin{cor}\label{cor:knr}
Let $n\geq 2$. Let 
\begin{equation*}
k= \begin{cases}
\frac{n-r}{2} & \text{if $n$ and $r$ are even}\\
n-r & \text{otherwise.}
\end{cases}
\end{equation*}
Then
\begin{equation*}
2l_{1}(n,r)= 
\begin{cases}
2\gcd(n,r) & \text{if $k$ is odd}\\
4\gcd(n,k) & \text{if $k$ is even.}
\end{cases}
\end{equation*}
\end{cor}

\begin{proof}
Suppose first that $n$ and $r$ are even, so $k=\frac{n-r}{2}$. If $k$ is odd then
$l_{1}(n,r)=\gcd(2r,n-r)=2\gcd \left( r,k \right)= 2\gcd \left( \frac{r}{2},\frac{n-r}{2}
\right)= \gcd \left( r,n-r \right)=\gcd(n,r)$. If $k$ is even then $l_{1}(n,r)= \gcd(2r,n-r)=
\gcd(2n,n-r)=2\gcd\left( n, \frac{n-r}{2}\right)= 2\gcd\left( n, k \right)$, which yields the
result in this case.

Now suppose that at least one of $n$ and $r$ is odd, so $k=n-r$. If $k$ is odd then
$l_{1}(n,r)=\gcd(2r,n-r)=\gcd(r,n-r)=\gcd(n,r)$. If $k$ is even then both $n$ and $r$ are odd, and
$l_{1}(n,r)=\gcd(2r,n-r)= 2\gcd(r,n-r)=2\gcd(n,n-r)=2 \gcd(n,k)$. 
\end{proof}

\section{Virtually cyclic groups}\label{sec:virtcyc}

We start this section by recalling the definition of a virtually cyclic group.
\begin{defn}
A group is \emph{virtually cyclic} if it contains a cyclic subgroup of finite index.
\end{defn}

\begin{rems}\mbox{}
\begin{enumerate}[(a)]
\item Every finite group is virtually cyclic.
\item Every infinite virtually cyclic group contains a normal subgroup of finite index.
\end{enumerate}
\end{rems}

The following criterion is due to Wall:
\begin{thm}[\cite{P,W}]\label{th:dd}
Let $G$ be a group. Then the following are equivalent.
\begin{enumerate}[(a)]
\item $G$ is a group with two ends.
\item $G$ is an infinite virtually cyclic group.
\item $G$ has a finite normal subgroup $F$ such that $G/F$ is $\Z$ or $\Z_{2} \ast \Z_{2}$.
\end{enumerate}
Equivalently, $G$ is of the form: 
\begin{enumerate}[(i)]
\item\label{it:semi} $F \rtimes \Z$, or 
\item\label{it:amalg} $G_{1} \ast_{F} G_{2}$, where $[G_{i} : F]=2$ for $i=1,2$.
\end{enumerate}
\end{thm}

If a virtually cyclic group $G$ satisifies~(i), we shall say that it is of \textit{Type~I}, while if
it satisifies~(ii), we shall say that it is of \textit{Type~II}.

We have already obtained the finite subgroups of $P_{n}(\rp)$ in \repr{finitepn}.
\reth{dd} allows us to establish a list of their \emph{possible} infinite virtually cyclic subgroups.
However, the difficulty is to decide whether the groups belonging to this list are effectively
realised as subgroups of $P_{n}(\rp)$. 

The following lemma gives a practical criterion for deciding whether a given infinite group is
virtually cyclic of Type~II, and shall be applied frequently in what follows.

\begin{lem}\label{lem:type2crit}
Let $G=G_{1} \ast_{F} G_{2}$ be a virtually cyclic group of Type~II, and let
$\map {\phi}{G_{1} \ast_{F} G_{2}}[H]$ be a homomorphism such that the restriction of $\phi$
to each $G_i$ is injective. Then $\phi$ is injective if and only if $\phi(G)$ is infinite.
\end{lem} 

\begin{proof} 
Clearly the given condition is necessary. So suppose that $\phi(G)$ is infinite. By \reth{dd}, we
have the short exact sequence
\begin{equation*}
1\to F\to G_1 \ast_F G_2 \to \Z_2 \ast \Z_2 \to 1.
\end{equation*}
Let $x$ (resp.\ $y$) denote the generator of the first (resp.\ second) copy of $\Z_{2}$. Then
$\ang{xy^{-1}}\cong \Z$ is a normal subgroup of $\Z_2 \ast \Z_2$, and $\Z_2 \ast \Z_2 \cong \Z
\rtimes \Z_{2}$, where $\Z_{2}$ may be taken as being generated by the $\ang{xy^{-1}}$-coset of $x$.

Let $K$ be the preimage of $\ang{xy^{-1}}$ under the projection $G \to \Z_2 \ast
\Z_2$. Then we have a short exact sequence
\begin{equation*}
1\to F\to K \to \Z \to 1,
\end{equation*}
and thus $K \cong F\rtimes \Z$ is virtually cyclic of Type~I, and is of index $2$ in $G$. This
gives rise to the following commutative diagram of short exact sequences:
\begin{equation*}
\xymatrix{%
& & 1 \ar[d] & 1 \ar[d] &\\
1\ar[r] & F \ar@{^{(}->}[r] \ar@{=}[d] & K \ar[r]
\ar@{^{(}->}[d] & \Z \ar[r] \ar@{^{(}->}[d] & 1\\
1\ar[r] & F \ar@{^{(}->}[r] & G \ar[r] \ar[d] &  \Z_{2} \ast \Z_{2} \ar[r] \ar[d]  &
1\\
& & \Z_{2} \ar@{=}[r] \ar[d] & \Z_{2} \ar[d] &\\
& & 1 & 1 &}
\end{equation*}

Let $w\in G_{1} \setminus F$, and consider the image of $w$ in $G_1 \ast_F G_2$, which by abuse of
notation we also denote by $w$. Since $w$ is of finite order it cannot belong to $K$ (for then it
would have to be mapped onto $0$ in $\Z$, which is impossible), and hence $G= K \coprod wK$. So by
hypothesis, $\phi(K)$ is infinite.

Notice that the result that we are aiming to prove is true for $K$ i.e.\ if we consider $\phi\lvert_{K}$,
and suppose that the restriction of $\phi$ to $F$ is injective then the fact that $\phi(K)$ is
infinite implies that $\phi\lvert_{K}$ is injective. Indeed, identifying $K$ with $F \rtimes \Z$, let
$k=(x,m)\in K$ belong to the kernel of $\phi\lvert_{K}$, where $x\in F$ and $m\in \Z$. Since
$\phi(K)$ is infinite, $\phi\lvert_{\Z}$ is injective. If $m\neq 0$ then $e_{H}= \phi(k)=
\phi((x,0)) \ldotp \phi((e_{G}, m))$. But $\phi((e_{G}, m))$ is of infinite order, while
$\phi((x,0))$ is of finite order, a contradiction. So $m=0$, $x=e_{F}$ by injectivity of
$\phi\lvert_{F}$, and so $\phi\lvert_{K}$ is injective.

Now $F \subset G_{1}$, so $\phi\lvert_{F}$ is injective. It follows then from the previous
paragraph that $\phi\lvert_{K}$ is injective. Furthermore, $\phi\lvert_{wK}$ is injective since
any two elements differ by an element of $K$.

Finally, to prove the result, it suffices to show that for all $k,k'\in K$, $\phi(k) \neq
\phi(wk')$, or equivalently that for all $k'' \in K$, $\phi(w)\neq \phi(k'')$. Suppose on the
contrary that there exists $k''\in K$ such that $\phi(w)=\phi(k'')$. Then $\phi(w^2)=\phi (k''^2)$.
Now $w^2\in K$ since $K$ is of index $2$ in $G$, and so $w^2=k''^2$ by injectivity of
$\phi\lvert_{K}$. Thus $k'' \in K$ is of finite order, and so belongs to $F$. Hence $w,k''\in
G_{1}$, $w\neq k''$ (as $w\in G_{1}\setminus F$), and so $\phi(w)\neq \phi (k'')$ by injectivity
of $\phi\lvert_{G_{1}}$. This yields a contradiction. Hence $\phi$ is injective.
\end{proof}

\section{Virtually cyclic subgroups of $P_{n}(\rp)$}

We now turn to the study of the virtually cyclic subgroups of $P_{n}(\rp)$. As $P_{1}(\rp) \cong
\Z_{2}$ and $P_{2}(\rp)\cong \quat$, it is trivial to determine their virtually cyclic subgroups. We
thus suppose from now on that $n\geq 3$. Since the structure of the finite subgroups differ for $n=3$
and $n\geq 4$, we treat these two cases separately. Further, by \repr{finitepn}, up to isomorphism,
we already know their finite subgroups. So in what follows, we shall seek their infinite virtually
cyclic subgroups.

\subsection{Virtually cyclic subgroups of $P_{3}(\rp)$}\label{sec:virtcycp3}

In this section, we prove \reth{vcp3} which is the case $n=3$ of \reth{vcp}. 
\begin{thm}\label{th:vcp3}
Up to isomorphism, the infinite virtually cyclic subgroups of $P_{3}(\rp)$ are $\Z$,
$\Z_{2}\times \Z$ and $\Z_{4} \ast_{\Z_{2}} \Z_{4}$.
\end{thm}

The key result needed to determine the Type~I subgroups of $P_{3}(\rp)$ is the following:
\begin{prop}\label{prop:p3noz4z}
The pure braid group $P_{3}(\rp)$ has no subgroup isomorphic to $\Z_{4}\times \Z$.
\end{prop}

\begin{proof}[Proof of \repr{p3noz4z}]
Suppose that $P_{3}(\rp)$ possesses a subgroup $G$ which is isomorphic to $\Z_{4} \times \Z$.
By~\cite{vB}, the Fadell-Neuwirth short exact sequence
\begin{equation*}
1 \to P_{1}(\rp \setminus \brak{x_{1},x_{2}}) \to P_{3}(\rp) \stackrel{p_{\ast}}{\to} P_{2}(\rp) \to 1
\end{equation*}
splits. Using the notation of page~765 of~\cite{GG2}, the kernel is a free group which we write as
$\F[2]=\F[2](x,y)$, where $x=\rho_{3}$ and $y=\rho_{3}^{-1} B_{2,3}$. The quotient $P_{2}(\rp)$ is
isomorphic to $\quat$, and is generated by $\rho_{1}$ and $\rho_{2}$ which are of order $4$. In
Corollary~11 of~\cite{GG2}, we exhibited an explicit section $s_{\ast}$ for $p_{\ast}$ given by
$s_{\ast}(\rho_{1})=\tau_{1}$, $s_{\ast}(\rho_{2})=\tau_{2}$, and
$s_{\ast}(\rho_{1}\rho_{2})=\tau_{3}=\tau_{1}\tau_{2}$, where the action on the kernel is as follows:
\begin{align*}
& \tau_{1} x \tau_{1}^{-1}= y & &  \tau_{2} x \tau_{2}^{-1}= y^{-1}  & & \tau_{3} x
\tau_{3}^{-1} = x^{-1}\\
& \tau_{1} y \tau_{1}^{-1}= x  & & \tau_{2} y \tau_{2}^{-1}= x^{-1}
& & \tau_{3} x \tau_{3}^{-1}= y^{-1}.
\end{align*}

Further, by Proposition~21 of~\cite{GG2}, there are precisely three conjugacy classes in
$P_{3}(\rp)$ of elements of order $4$ whose representatives may be taken to be $\tau_{1},\tau_{2}$
and $\tau_{3}$. So conjugating $G$ by an element of $P_{3}(\rp)$ if necessary, and identifying
$P_{3}(\rp)$ with $\F[2]\rtimes \quat$, we may suppose that its $\Z_{4}$ factor is generated by 
$(1,\tau_{i})$ for some $i=1,2,3$. Let $(z,\xi)$ generate the $\Z$-factor of $G$ for some $z\in
\F[2]$ and $\xi\in \quat$. Since $G$ is the direct product of these two factors, we have:
\begin{align*}
(1,\tau_{i})\ldotp (z,\xi) &= (z,\xi)\ldotp (1,\tau_{i})\\
(\tau_{i}z \tau_{i}^{-1}, \tau_{i}\xi) &= (z,\xi \tau_{i}).
\end{align*}
Thus $z=\tau_{i}z \tau_{i}^{-1}$. But from the form of the action of the $\tau_{i}$ on $x$ and $y$
given above, this is not possible unless $z=1$, in which case $(z,\xi)$ is of finite order, a
contradiction.
\end{proof}

\begin{proof}[Proof of \reth{vcp3}.]
Let $G$ be an infinite virtually cyclic subgroup of $P_{3}(\rp)$. We first suppose that it is of
Type~I. By \reth{dd}, $G$ must be isomorphic to one of the following groups: $\Z$, $\Z_{2}\rtimes
\Z$, $\Z_{4} \rtimes \Z$ or $\quat \rtimes \Z$. Clearly the two groups $\Z$ and $\Z_{2}\times
\Z$ are realised as subgroups of $P_{3}(\rp)$, and there is no non-trivial semi-direct product
$\Z_{2}\rtimes \Z$. We saw in \repr{p3noz4z} that $G$ cannot be isomorphic to $\Z_{4} \times \Z$. 
Since $\Z_{4}$ is a subgroup of $\quat$, it follows that $\quat \times \Z$ cannot be realised as a
subgroup of $P_{3}(\rp)$ either.

The possible (non-trivial) semi-direct products are not realised either as subgroups of $P_{3}(\rp)$.
Indeed, since $\aut{\Z_{4}}\cong \Z_{2}$, it follows that the subgroup $\Z_{4}\rtimes 2\Z$ would in
fact be a direct product, abstractly isomorphic to $\Z_{4} \times \Z$. Similarly, since
$\aut{\quat}\cong \sn[4]$, the subgroup $\quat \rtimes 12\Z$ would be a direct product, abstractly
isomorphic to $\quat \times \Z$. Thus $\Z$ and $\Z_{2}\times \Z$ are the only groups realised as
Type~I subgroups of $P_{3}(\rp)$.

Now suppose that $G$ is of Type~II. By \reth{dd}, $G$ is isomorphic to one of $\Z_{2} \ast
\Z_{2}$, $\Z_{4} \ast_{\Z_{2}} \Z_{4}$ or $\quat \ast_{\Z_{4}} \quat$. The first case is ruled out
since $P_{3}(\rp)$ has just one element of order $2$.  If $P_{3}(\rp)$ had a subgroup isomorphic to
$\quat \ast_{\Z_{4}} \quat$ then by the proof of \relem{type2crit}, it would also have a subgroup
isomorphic to $\Z_{4}\rtimes \Z$ which as we just saw cannot happen. The realisation of $\Z_{4}
\ast_{\Z_{2}} \Z_{4}$ as a subgroup of $P_{n}(\rp)$ occurs for all $n\geq 3$, and will be proved below
in \repr{z4astz4}. This completes the proof of \reth{vcp3}.
\end{proof}

\subsection{Virtually cyclic subgroups of $P_{n}(\rp)$, $n\geq 4$}\label{sec:virtcycp4}

In this section, we classify the virtually cyclic subgroups of $P_{n}(\rp)$ in the general case. In
fact we obtain exactly the same result as in the case $n=3$.
\begin{thm}\label{th:vcp4}
Let $n\geq 4$. Up to isomorphism, the infinite virtually cyclic subgroups of $P_{n}(\rp)$ are $\Z$,
$\Z_{2}\times \Z$ and $\Z_{4} \ast_{\Z_{2}} \Z_{4}$.
\end{thm}

In order to prove \reth{vcp4}, we first state and prove the following two propositions concerning
the realisation of $\Z_{4} \ast_{\Z_{2}} \Z_{4}$ and the non-realisation of $\Z_{4}\times \Z$.

\begin{prop}\label{prop:z4astz4}
Let $n\geq 3$. Then $P_{n}(\rp)$ possesses a subgroup isomorphic to $\Z_{4} \ast_{\Z_{2}} \Z_{4}$.
\end{prop}

\begin{proof}[Proof of \repr{z4astz4}]
Let $H=\ang{\alpha,\beta}< P_{n}(\rp)$, where $\alpha,\beta$ are as defined in
\req{defalpha}. Since $\alpha,\beta$ are of order $4$ and $\alpha^2=\beta^2$, there is a natural
surjective homomorphism $\Z_{4} \ast_{\Z_{2}} \Z_{4} \to H$ that is injective on the
$\Z_{4}$-factors. Further $\rho_{n}=\alpha \beta^{-1} \in H$ is of infinite order, so $\ang{\alpha}
\cap \ang{\beta}= \ang{\ft}$ and $H$ is of infinite order. We conclude from \relem{type2crit} that
$H\cong \Z_{4} \ast_{\Z_{2}} \Z_{4}$.
\end{proof}

\begin{rem}\label{rem:z4z2z4}
A similar argument may be applied in order to show that if $n\geq 4$ then any two distinct subgroups
of $P_{n}(\rp)$ of order $4$ generate a subgroup isomorphic to $\Z_{4} \ast_{\Z_{2}} \Z_{4}$. The
fact that this subgroup is infinite follows directly from \repr{finitepn}.
\end{rem}

\begin{prop}\label{prop:noz4z}
For all $n\geq 1$, $P_{n}(\rp)$ has no subgroup isomorphic to $\Z_{4}\times \Z$.
\end{prop}

\begin{proof}
Since $P_{1}(\rp)$ and $P_{2}(\rp)$ are finite, we may suppose that $n\geq 3$. We argue by
induction on $n$. If $n=3$ then the result was proved previously in \repr{p3noz4z}. So suppose that
the result is true for $P_{n}(\rp)$, and let us consider $P_{n+1}(\rp)$. Suppose that $P_{n+1}(\rp)$
has a subgroup $H$ isomorphic to $\Z_{4} \times \Z$. Let $x,g\in P_{n+1}(\rp)$ generate respectively
the $\Z_{4}$- and $\Z$-factors. Since $P_{n+1}(\rp)$ is normal in $B_{n+1}(\rp)$, by conjugating $H$
by an element of $B_{n+1}(\rp)$ if necessary, we may suppose further by \repr{order4} that $x$ is
equal to one of $\alpha=a^{n+1}$ or $\beta=b^n$.

Consider the usual projection $\map{p_{\ast}}{P_{n+1}(\rp)}[P_{n}(\rp)]$, consisting in forgetting
the last string. Since $\ker{p_{\ast}}=P_{1}(\rp \setminus \brak{x_{1}, \ldots, x_{n}})$ is torsion
free, $p_{\ast}$ is injective on finite subgroups of $P_{n+1}(\rp)$, and so $x'=p_{\ast}(x)$ is of
order $4$. Let $g'=p_{\ast}(g)$. If $g'$ is of infinite order then $\ang{x',g'}$ is a subgroup of
$P_{n}(\rp)$ isomorphic to $\Z_{4} \times \Z$, which is ruled out by the induction hypothesis. So
$g'$ must be of finite order. Since $g'$ and $x'$ commute, they generate a finite Abelian subgroup of
$P_{n}(\rp)$ which by \repr{finitepn} must be isomorphic to $\Z_{4}$. In particular, $g'\in
\ang{x'}$, and there exists $0\leq j\leq 3$ such that $z=gx^{-j}\in \ker{p_{\ast}}$. Then $z$ is of
infinite order, $H$ is generated by $z$ and $x$, and is isomorphic to the direct product
$\ang{x}\times\ang{z}\cong\Z_{4}\times \Z$. So replacing $g$ by $z$ if necessary, we may suppose
that $g\in \ker{p_{\ast}}$. 

If such a $g$ were indeed to exist, it would be a non-trivial fixed point of the action of
conjugation by $x$ on $\ker{p_{\ast}}$. So let us study this action, and show that there are no such
fixed points, which will thus yield a contradiction. 

Now $\ker{p_{\ast}}= \pi_{1}(\rp \setminus \brak{x_{1},\ldots, x_{n}}, x_{n+1})$ is a free
group of rank~$n$ for which a basis is:
\begin{equation}\label{eq:basis}
\brak{\rho_{n+1}, B_{1,n+1}, \ldots, B_{n-1,n+1}},
\end{equation}
where for $i=1,\ldots,n$,
\begin{equation*}
B_{i,n+1}=\sigma_{n}\cdots \sigma_{i+1} \sigma_{i}^2 \sigma_{i+1}^{-1} \cdots
\sigma_{n}^{-1}.
\end{equation*}
By equation~(7) of~\cite{GG2}, we have
\begin{equation*}
\rho_{n+1}^{-2}=\sigma_{n} \cdots \sigma_{2} \sigma_{1}^2 \sigma_{2} \cdots \sigma_{n}, 
\end{equation*}
which is equal to $B_{1,n+1} \cdots B_{n-1,n+1} B_{n,n+1}$. So 
\begin{equation*}
\sigma_{n}^2= B_{n,n+1}= B_{n-1,n+1}^{-1} \cdots B_{1,n+1}^{-1}\rho_{n+1}^{-2}.
\end{equation*}
Further, by using the Artin relations of $B_{n}(\rp)$ (those involving just the $\sigma_{i}$), for
all $1\leq i\leq n-1$, we obtain the following useful identities:
\begin{equation}\label{eq:sigbin}
\left\{ \begin{aligned}
\sigma_{i} B_{i,n+1} \sigma_{i}^{-1} &= B_{i+1,n+1}\\
\sigma_{i}^{-1} B_{i,n+1} \sigma_{i} &= B_{i,n+1} B_{i+1,n+1} B_{i,n+1}^{-1}.
\end{aligned} \right.
\end{equation}

We now calculate the action of $\alpha$ and $\beta$ on each of the elements of the basis given in
\req{basis}. It was proved on page~777 of~\cite{GG2} that conjugation by $a^{-1}$ permutes cyclically
the $2(n+1)$~elements
\begin{equation*}
\sigma_{1}, \sigma_{2}, \ldots, \sigma_{n}, a^{-1} \sigma_{n}a, \sigma_{1}^{-1}, \sigma_{2}^{-1},
\sigma_{n}^{-1}, a^{-1} \sigma_{n}^{-1}a, 
\end{equation*}
and permutes cyclically the $2(n+1)$~elements
\begin{equation*}
\rho_{1}, \rho_{2}, \ldots, \rho_{n}, \rho_{n+1}, \rho_{1}^{-1}, \rho_{2}^{-1}, \ldots, \rho_{n}^{-1},
\rho_{n+1}^{-1}.
\end{equation*}
In particular, conjugation by $\alpha^{-1}=a^{-(n+1)}$ sends each of $\sigma_{i}$, $1\leq i\leq n$,
and $\rho_{j}$, $1\leq j\leq n+1$, to its inverse. By induction, we have that for all $1\leq i\leq n-1$, 
\begin{equation*}
\alpha^{-1} B_{i,n+1} \alpha = \rho_{n+1}^2 B_{1,n+1} \cdots B_{i-1,n+1} B_{i,n+1}^{-1}
B_{i-1,n+1}^{-1} \cdots B_{1,n+1}^{-1} \rho_{n+1}^{-2}.
\end{equation*}
To see this, first let $i=1$. Then:
\begin{align*}
\alpha^{-1} B_{1,n+1} \alpha &= \alpha^{-1} \sigma_{n}\cdots \sigma_{2}\sigma_{1}^2
\sigma_{2}^{-1} \cdots \sigma_{n}^{-1} \alpha\\
&= \sigma_{n}^{-1}\cdots \sigma_{2}^{-1} \sigma_{1}^{-2} \sigma_{2}^{-1} \cdots \sigma_{n}^{-1}
\ldotp \sigma_{n} \cdots \sigma_{2} \sigma_{1}^{-2} \sigma_{2}^{-1} \cdots \sigma_{n}^{-1} 
\ldotp \sigma_{n}\cdots \sigma_{2}\sigma_{1}^2 \sigma_{2} \cdots \sigma_{n}\\
&= \rho_{n+1}^2 B_{1,n+1}^{-1} \rho_{n+1}^{-2},
\end{align*}
using equation~(7) of~\cite{GG2}. Now suppose by induction that the result holds for $1\leq i\leq n-2$.
Then using \req{sigbin}, we have
\begin{align*}
\alpha^{-1} B_{i+1,n+1} \alpha &= \alpha^{-1} \sigma_{i} B_{i,n+1} \sigma_{i}^{-1} \alpha\\
&= \sigma_{i}^{-1}\left( \rho_{n+1}^2 B_{1,n+1} \cdots B_{i-1,n+1} B_{i,n+1}^{-1}
B_{i-1,n+1}^{-1} \cdots B_{1,n+1}^{-1} \rho_{n+1}^{-2} \right) \sigma_{i}\\
&=  \rho_{n+1}^2 B_{1,n+1} \cdots B_{i-1,n+1} \sigma_{i}^{-1} B_{i,n+1}^{-1} \sigma_{i}
B_{i-1,n+1}^{-1} \cdots B_{1,n+1}^{-1} \rho_{n+1}^{-2}\\
&=  \rho_{n+1}^2 B_{1,n+1} \cdots B_{i-1,n+1}  B_{i,n+1} B_{i+1,n+1}^{-1} B_{i,n+1}^{-1}
B_{i-1,n+1}^{-1} \cdots B_{1,n+1}^{-1} \rho_{n+1}^{-2},
\end{align*}
as required.

In other words, the action of conjugating each of these generators by $\alpha^{-1}$
is equal to $\iota_{\rho_{n+1}^2}\circ \phi$, where $\iota_{\rho_{n+1}^2}$ is conjugation by
$\rho_{n+1}^2$, and $\phi$ is the automorphism of $\ker{p_{\ast}}$ given by:
\begin{align*}
\phi(\rho_{n+1}) &= \rho_{n+1}^{-1}\\
\phi(B_{i,n+1}) &= B_{1,n+1} \cdots B_{i-1,n+1} B_{i,n+1}^{-1}
B_{i-1,n+1}^{-1} \cdots B_{1,n+1}^{-1},
\end{align*}
for $i=1,\ldots, n-1$. Notice further that $\phi$ induces an automorphism of the
subgroup $K=\ang{B_{1,n+1}, \ldots, B_{n-1,n+1}}$ of $\ker{p_{\ast}}$.

Since $\alpha= \rho_{n+1} \beta$, for all $y\in \ker{p_{\ast}}$,
\begin{equation*}
\beta^{-1} y \beta= \alpha^{-1} \rho_{n+1} y \rho_{n+1}^{-1} \alpha= \rho_{n+1}^{-1} \alpha^{-1} y
\alpha \rho_{n+1},
\end{equation*}
and so conjugation by $\beta^{-1}$ on each of the given basis elements is equal to
$\iota_{\rho_{n+1}}\circ \phi$.

Suppose first that $x=\alpha$, and let $z=\rho_{n+1}^{\epsilon_{0}} w_{1} \rho_{n+1}^{\epsilon_{1}}
\cdots w_{k-1} \rho_{n+1}^{\epsilon_{k}}$ be a fixed point of
the action of conjugation by $x$, where $k\geq 0$, for $i=1,\ldots, k-1$,
$w_{i}\in K$ is a non-trivial word, $\epsilon_{j}\in \Z$ for $j=0,1,\ldots, k$,
and $\epsilon_{j}\neq 0$ for $j=1,\ldots,k-1$. Now $\alpha^{-1} \rho_{n+1} \alpha =
\rho_{n+1}^{-1}$, so we must have $k\geq 1$, and then
\begin{equation*}
\alpha^{-1}z \alpha= \rho_{n+1}^{2-\epsilon_{0}} \phi(w_{1}) \rho_{n+1}^{-\epsilon_{1}} \cdots
\phi(w_{k-1}) \rho_{n+1}^{-\epsilon_{k}-2}.
\end{equation*}
Since $\phi$ induces an automorphism of $K$, there is no cancellation between successive terms
of this expression. It thus follows that $k=1$, $\epsilon_{0}=1$ and $\epsilon_{1}=-1$ and
$\phi(w_{1})=w_{1}$. Now $\phi\lvert_{K}=\iota_{B_{1,n+1}}\circ \phi'$, where $\phi'$ is the
automorphism of $K$ defined by:
\begin{align*}
\phi'(B_{1,n+1}) &= B_{1,n+1}^{-1}\\
\phi'(B_{i,n+1}) &= B_{2,n+1} \cdots B_{i-1,n+1} B_{i,n+1}^{-1} B_{i-1,n+1}^{-1} \cdots B_{2,n+1}^{-1},
\end{align*}
for $i=2,\ldots, n-1$. Similarly, $\phi'$ restricts to an automorphism of the subgroup
$K'=\ang{B_{2,n+1}, \ldots, B_{n-1,n+1}}$.

Set $w_{1}=B_{1,n+1}^{\delta_{0}} v_{1} B_{1,n+1}^{\delta_{1}}
\cdots v_{t-1} B_{1,n+1}^{\delta_{t}}$, where $t\geq 0$, for $i=1,\ldots, t-1$,
$v_{i}\in K'$ is a non-trivial word, $\delta_{j}\in \Z$ for $j=0,1,\ldots, t$,
and $\delta_{j}\neq 0$ for $j=1,\ldots,t-1$. Now $\phi(B_{1,n+1}) =
B_{1,n+1}^{-1}$, so we must have $t\geq 1$, and then
\begin{equation*}
\phi(w_{1})= B_{1,n+1}^{1-\delta_{0}} \phi'(v_{1}) B_{1,n+1}^{-\delta_{1}}
\cdots \phi'(v_{t-1}) B_{1,n+1}^{-\delta_{t}-1}.
\end{equation*}
Again, since $K'$ is invariant under $\phi'$, there is no cancellation between successive terms in
this expression. But $w_{1}$ is a (non-trivial) fixed point of $\phi$, so we must have
$2\delta_{0}=1$ which yields a contradiction.  Hence there exists no $g\in P_{n+1}(\rp)$ of infinite
order which commutes with $\alpha$.
We now repeat the same analysis for $\beta$, from which we obtain:
\begin{equation*}
\beta^{-1}z \beta= \rho_{n+1}^{1-\epsilon_{0}} \phi(w_{1}) \rho_{n+1}^{-\epsilon_{1}} \cdots
\phi(w_{k-1}) \rho_{n+1}^{-\epsilon_{k}-1}.
\end{equation*}
Thus $\epsilon_{0}\notin \Z$, a contradiction. Hence there exists no $g\in P_{n+1}(\rp)$
of infinite order that commutes with $\beta$. 
\end{proof}

\begin{proof}[Proof of \reth{vcp4}.]
Let $H$ be an infinite virtually cyclic subgroup of $P_{n}(\rp)$. If $H$ is of Type~I then by
\repr{finitepn}, $H$ is isomorphic to either $\Z$, $\Z_{2} \times \Z$ (both of which are realised),
$\Z_{4}\times \Z$ or $\Z_{4} \rtimes \Z$. But $\Z_{4}\times \Z$ and $\Z_{4} \rtimes \Z$ are ruled out
by \repr{noz4z}. If $H$ is of Type~II then either $H\cong \Z_{2} \ast \Z_{2}$, which is not possible
since $P_{n}(\rp)$ has only one subgroup isomorphic to $\Z_{2}$, or $H\cong \Z_{4} \ast_{\Z_{2}}
\Z_{4}$, which as we saw in \repr{z4astz4}, is realised. This completes the proof of \reth{vcp4}, as
well as that of \reth{vcp}.
\end{proof}

\begin{rem}
We saw that in $B_{n}(\rp)$, there are two conjugacy classes of subgroups of order $4$ lying in
$P_{n}(\rp)$. It would be interesting to know how many such conjugacy classes exist in $P_{n}(\rp)$.
\end{rem}

As a corollary of \reth{vcp}, we obtain the following result:
\begin{cor}
Let $n\geq 4$. Then the virtually cyclic subgroups of $P_{n}(\rp)$ are $\brak{e}$, $\Z_{2}$,
$\Z_{4}$, $\Z$, $\Z_{2} \times \Z$ and $\Z_{4} \ast_{\Z_{2}} \Z_{4}$.\qed
\end{cor}

We also obtain \reco{centralx}:
\begin{proof}[Proof of \reco{centralx}]
Clearly $\ang{x}\subseteq Z_{P_{n}(\rp)}(x)$. So suppose that $z\in Z_{P_{n}(\rp)}(x)$. It follows
from \repr{noz4z} that $z$ is of finite order. Then $x$ and $z$ generate a finite Abelian subgroup
of $P_{n}(\rp)$ which contains $\ang{x}$. So $\ang{x,z}$ is isomorphic to $\Z_{4}$, and thus $z\in
\ang{x}$.
\end{proof}

\begin{rem}
For $x \in P_{n}(\rp)$ of order $4$, we believe that the centraliser of $Z_{B_{n}(\rp)}(x)$ in
$B_{n}(\rp)$ is finite cyclic of order $4n$ or $4(n-1)$ depending on whether $x$ is conjugate in
$B_{n}(\rp)$ to $\alpha$ or $\beta$. In order to justify this, notice first that $Z_{B_{n}(\rp)}(x)$
cannot contain any element $z$ of infinite order, since there would exist some $1\leq j\leq n!$ such
that $z^j\in P_{n}(\rp)$, and so $z^j\in  Z_{P_{n}(\rp)}(x)$, a contradiction. Then
$Z_{B_{n}(\rp)}(x)$ contains only elements of finite order, and so is itself finite. To prove this,
suppose that $Z_{B_{n}(\rp)}(x)$ is infinite. Then $K=Z_{B_{n}(\rp)}(x)\cap P_{n}(\rp)$ is infinite. 
\begin{enumerate}[--]
\item If $n\geq 4$ then $K$ contains two distinct copies of $\Z_{4}$, and by \rerem{z4z2z4} they
together generate a copy of $\Z_{4} \ast_{\Z_{2}} \Z_{4}$ which contains elements of infinite order.
\item If $n=3$ then $K$ contains an infinite number of elements of order $4$. Identifying $P_{3}(\rp)$ 
with $\F[2] \rtimes \quat$ as in the proof of \repr{p3noz4z}, there exist $i\in
\brak{1,2,3}$ and $z_{1},z_{2}\in \F[2](x,y)$, $z_{1}\neq z_{2}$, such that $(z_{1}, \tau_{i}),
(z_{2}, \tau_{i}) \in K$. But then 
\begin{equation*}
(z_{1}, \tau_{i}) \ldotp (z_{2}, \tau_{i})^{-1}= (z_{1}, \tau_{i}) \ldotp (\tau_{i}^{-1}z_{2}^{-1}
\tau_{i}, \tau_{i}^{-1})= (z_{1}z_{2}^{-1}, 1)\in K.
\end{equation*}
\end{enumerate}
In both cases, we deduce a contradiction, and thus $Z_{B_{n}(\rp)}(x)$ is finite. Further, by
conjugating $x$ by an element of $B_{n}(\rp)$ if necessary, we may suppose that $x=\alpha=a^n$ or
$x=\beta=b^{n-1}$, in which case $\ang{a}\subseteq Z_{B_{n}(\rp)}(x)$ or $\ang{b}
\subseteq Z_{B_{n}(\rp)}(x)$ respectively. Investigations into the structure of the finite
subgroups of $B_{n}(\rp)$ suggest that $\ang{a}$ and $\ang{b}$ are maximal finite subgroups of
$B_{n}(\rp)$. If this is indeed the case then $Z_{B_{n}(\rp)}(x)=\ang{a}$ or
$Z_{B_{n}(\rp)}(x)=\ang{b}$ as claimed. A similar argument would show that the normaliser of
$\ang{x}$ in $B_{n}(\rp)$ is also finite cyclic of order $4n$ or $4(n-1)$.
\end{rem}

\end{document}